\documentclass[11pt]{article}
\usepackage{amssymb,amsthm,amsmath}
\def\qed{\hfill $\Box$}

\newcommand\pf{\smallbreak\noindent \texttt{Proof}. }

\begin{document}

\newtheorem{thm}{Theorem}[section]
\newtheorem{prop}[thm]{Proposition}
\newtheorem{lem}[thm]{Lemma}
\newtheorem{cor}[thm]{Corollary}
\newtheorem{ex}[thm]{Example}
\renewcommand{\thefootnote}{*}

\title{\bf On Leibniz algebras, whose subalgebras are either ideals or self-idealizing}

\author{\textbf{L.A.~Kurdachenko, A.A.~Pypka}\\
Department of Geometry and Algebra\\
Oles Honchar Dnipro National University, Dnipro, Ukraine\\
{\small e-mail: lkurdachenko@gmail.com, sasha.pypka@gmail.com}\\
\textbf{I.Ya.~Subbotin}\\
Department of Mathematics and Natural Sciences\\
National University, Los Angeles, USA\\
{\small e-mail: isubboti@nu.edu}}
\date{}

\maketitle

\begin{abstract}
A subalgebra $S$ of a Leibniz algebra $L$ is called self-idealizing in $L$ if it coincides with its idealizer $\mathrm{I}_{L}(S)$. In this paper we study the structure of Leibniz algebras, whose subalgebras are either ideals or self-idealizing.
\end{abstract}



\thispagestyle{empty}

\section{Introduction}
Let $L$ be an algebra over a field $F$ with the binary operations $+$ and $[\ ,\ ]$. Then $L$ is called a \textit{left Leibniz algebra} if it satisfies the left Leibniz identity
$$[[a,b],c]=[a,[b,c]]-[b,[a,c]]$$
for all $a,b,c\in L$. We will also use another form of this identity:
$$[a,[b,c]]=[[a,b],c]+[b,[a,c]].$$

Leibniz algebras appeared first in the paper of A.~Blokh \cite{BA1965}, in which he called them the \textit{$D$-algebras}. However, in that time, these works were not in demand, and they have not been properly developed. A real interest to Leibniz algebras appears only after two decades. The term ``Leibniz algebra'' appears in the book of J.-L.~Loday \cite{LJ1992}, and the article of J.-L.~Loday \cite{LJ1993}. In \cite{LP1993} J.-L.~Loday and T.~Pirashvili began a study of the properties of Leibniz algebras. The theory of Leibniz algebras has developed very intensively in many different directions. Some of the results of this theory were presented in the book \cite{AOR2020}. Note that Lie algebras are a partial case of Leibniz algebras. Conversely, if $L$ is a Leibniz algebra, in which $[a,a]= 0$ for every element $a\in L$, then it is a Lie algebra. Thus, Lie algebras can be characterized as anticommutative Leibniz algebras. In this regard, a parallel with associative structures, such as, for example, groups and associative rings, comes to mind. There we can see a significant difference between Abelian and non-Abelian groups, between commutative and non-commutative rings. The differences between Lie algebras and Leibniz algebras became visible when one considers the first natural types of Leibniz algebras. For example, cyclic Lie algebras has dimension $1$, however the structure of cyclic Leibniz algebras is much more sophisticated. The structure of cyclic Leibniz algebras has been described in \cite{CKS2017}. Another example of this sort: Lie algebras, every subalgebra of which is ideal, is abelian, while in the case of Leibniz algebras, we have a different more complex case. Leibniz algebras, whose subalgebras are ideals, have been described in \cite{KSS2017}. We face the same situation when studying other types of Leibniz algebras (see papers \cite{KOS2019,KSS2018A,KSY2018,KSY2020}).

In the current paper, the term ``Leibniz algebra'' stands for a left Leibniz algebra that is not a Lie algebra.

Let $L$ be a Leibniz algebra over a field $F$, $A$ and $K$ be subalgebras of $L$. The \textit{left idealizer} of $A$ in $K$ is defined by the following rule:
$$\mathrm{I}_{K}^{left}(A)=\{x\in K|\ [x,a]\in A\ \mbox{for all elements }a\in A\}.$$

The left idealizer of $A$ in $K$ is a subalgebra of $K$. Indeed, let $x,y\in\mathrm{I}_{K}^{left}(A)$, $a\in A$, $\alpha\in F$, then
\begin{equation*}
\begin{split}
&[x-y,a]=[x,a]-[y,a]\in A,\ [\alpha x,a]=\alpha[x,a]\in A,\\
&[[x,y],a]=[x,[y,a]]-[y,[x,a]]\in A.
\end{split}
\end{equation*}

Similarly the \textit{right idealizer} of $A$ in $K$ is defined by the following rule:
$$\mathrm{I}_{K}^{right}(A)=\{x\in K|\ [a,x]\in A\ \mbox{for all elements }a\in A\}.$$

Unlike left idealizer, the right idealized of $A$ in $K$ does need not to be a subalgebra. This is shown in the corresponding example from \cite[Example~1.7]{BD2011}.

The \textit{idealizer} of $A$ in $K$ is defined by the following rule:
$$\mathrm{I}_{K}(A)=\{x\in K|\ [x,a],[a,x]\in A\ \mbox{for all elements }a\in A\}=\mathrm{I}_{K}^{left}(A)\cap\mathrm{I}_{K}^{right}(A).$$

The idealizer of $A$ in $K$ is a subalgebra of $K$. Indeed, let $x,y\in\mathrm{I}_{K}(A)$, $a\in A$, $\alpha\in F$. As we did it above, we can prove that $x-y,\alpha x\in\mathrm{I}_{K}(A)$. Further,
\begin{equation*}
\begin{split}
&[[x,y],a]=[x,[y,a]]-[y,[x,a]]\in A,\\
&[a,[x,y]]=[[a,x],y]+[x,[a,y]]\in A.
\end{split}
\end{equation*}

For arbitrary subalgebra $A$ of a Leibniz algebra $L$ we have the following ascending series
$$A=A_{0}\leqslant A_{1} \leqslant\ldots A_{\alpha}\leqslant A_{\alpha+ 1} \leqslant\ldots A_{\gamma}$$
where $A_{1} =\mathrm{I}_{L}(A)$, $A_{\alpha+ 1} =\mathrm{I}_{L}(A_{\alpha})$ for all ordinals $\alpha$, $A_{\lambda}=\bigcup_{\mu<\lambda}A_{\mu}$ for all limit ordinals $\lambda$ and $A_{\gamma}=\mathrm{I}_{L}(A_{\gamma})$. This series is called the \textit{upper idealized series} of $A$ in $L$. If $\gamma= 1$, then we obtain the following two cases:
$$A_{1} =\mathrm{I}_{L}(A)=L\ \mbox{or }A_{1} =\mathrm{I}_{L}(A)=A.$$
The following two types of subalgebras correspond to these two cases: $A$ is an ideal of $L$ or $A$ is self-idealizing in $L$ (that is $A=\mathrm{I}_{L}(A))$. Therefore, the following natural question appears:
\begin{center}
\textit{What is a Leibniz algebra, every subalgebra of which either is an ideal or self-idealizing?}
\end{center}

The purpose of this work is to develop a precise description of such Leibniz algebras.

The first type of such algebras are the Leibniz algebras, whose subalgebras are ideals. These algebras have been described in \cite{KSS2017}.

A Leibniz algebra $L$ is called an \textit{extraspecial}, if it satisfies the following condition:
\begin{itemize}
\item the center $\zeta(L)$ of $L$ is non-trivial and of dimension $1$,
\item $L/\zeta(L)$ is abelian.
\end{itemize}

An extraspecial algebra $E$ is called a \textit{strong extraspecial} algebra, if $[x,x]\neq 0$ for each element $x\not\in\zeta(E)$.

Leibniz algebra $L$, whose subalgebras are ideals, has the following structure: $L=E\oplus Z$ where $Z$ is a subalgebra of the center of $L$ and $E$ is a strong extraspecial algebra.

The study of Leibniz algebras, whose subalgebras are either ideal or self-idealizing, consist of the following steps.

Let $L$ be a Leibniz algebra over a field $F$. Then $L$ includes the greatest locally nilpotent ideal \cite[Corollary~C1]{KSS2018}. This ideal is called the \textit{locally nilpotent radical} of $L$ and will be denoted by $\mathrm{Ln}(L)$.

At the first stage, we study Leibniz algebras, whose locally nilpotent radical is abelian and not cyclic.

\medskip

\noindent
{\bf Theorem~A}\ {\it Let $L$ be a Leibniz algebra over a field $F$, whose subalgebras are either ideals or self-idealizing. Suppose that the locally nilpotent radical $A$ of $L$ is abelian. If $A$ is not cyclic, then the following conditions hold:
\begin{enumerate}
\item[\upshape(i)] $A=\zeta^{left}(L)$;
\item[\upshape(ii)] $L=A\oplus W$ where $W$ is a subalgebra of dimension $1$, $W=Fw$, $W=\mathrm{I}_{L}(W)$;
\item[\upshape(iii)] there exists a non-zero element $\sigma\in F$ such that $[w,a]=\sigma a$ for every element $a\in A$.
\end{enumerate}

Conversely, if a Leibniz algebra $L$ satisfies these conditions, then every subalgebra of $L$ is an ideal of $L$ or self-idealizing in $L$.}

At the next stage, it is natural to study the case when the locally nilpotent radical is non-abelian and not cyclic. Here we have the following results.

\medskip

\noindent
{\bf Theorem~B1}\ {\it Let $L$ be a Leibniz algebra over a field $F$, whose subalgebras are either ideals or self-idealizing. Suppose that $L\neq\mathrm{Ln}(L)$, $\mathrm{Ln}(L)$ is a non-cyclic subalgebra and $\mathrm{Ln}(L)$ is not abelian. If $[\mathrm{Ln}(L),\mathrm{Ln}(L)]\leqslant\zeta(L)$, then $\mathrm{char}(F)= 2$ and the following conditions hold:
\begin{enumerate}
\item[\upshape(i)] $K=\mathrm{Ln}(L)$ is a strong extraspecial subalgebra, moreover $[x,x]\neq 0$ for each element $x\not\in\mathrm{Leib}(L)$;
\item[\upshape(ii)] $[K,K]=\zeta(L)=\mathrm{Leib}(L)$;
\item[\upshape(iii)] $L=K+\langle v\rangle$ where $[v,v]=\eta z$ for some non-zero element $\eta\in F$, $[v+\zeta(K),x+\zeta(K)]=x+\zeta(K)$ for each element $x\in K\setminus\zeta(K)$.
\end{enumerate}}

\noindent
{\bf Theorem~B2}\ {\it Let $L$ be a Leibniz algebra over a field $F$, whose subalgebras are either ideals or self-idealizing. Suppose that $L\neq\mathrm{Ln}(L)=K$ is not cyclic and not abelian. If $\zeta(L)$ does not include $[K,K]$, then the following conditions hold:
\begin{enumerate}
\item[\upshape(i)] $\mathrm{char}(F)\neq 2$;
\item[\upshape(ii)] every subalgebra of $K$ is an ideal of $L$;
\item[\upshape(iii)] $\mathrm{Ln}(L)$ is a strong extraspecial subalgebra;
\item[\upshape(iv)] $[K,K]=\zeta(K)=\mathrm{Leib}(L)=\langle z\rangle$ has dimension $1$;
\item[\upshape(v)] $L=K+\langle v\rangle$ where $[v+\zeta(K),x+\zeta(K)]=x+\zeta(K)$ for each element $x\in K\setminus\zeta(K)$, $[v,z]= 2z$ and $[v,v]=\nu z$ for some element $\nu\in F$ ($\nu$ can be zero).
\end{enumerate}}

The last step is to study the case when the locally nilpotent radical is cyclic. In this case, its dimension is $1$ or $2$. The following theorem describes the second case.

\medskip

\noindent
{\bf Theorem~C}\ {\it Let $L$ be a Leibniz algebra over a field $F$, whose subalgebras are either ideals or self-idealizing. Suppose that $\mathrm{Ln}(L)$ is a cyclic subalgebra, having dimension $2$. Then either $L\neq\mathrm{Ln}(L)$, or $\mathrm{char}(F)= 2$ and $L$ has a basis $\{z,a,v\}$ such that $[z,z]=[z,a]=[a,z]=[z,v]=[v,z]= 0$, $[a,a]=z$, $[v,v]=\eta z$, $[v,a]=a+\lambda z$, $[a,v]=a+\mu z$, $\eta,\lambda,\mu\in F$ and a polynomial $X^{2}+(\mu+\lambda)X+\eta$ has no root in $F$.}

The case when the dimension of the locally nilpotent radical is $1$ reduces to the study of Lie algebras, whose abelian subalgebras have dimension $1$. This situation requires a separate consideration. In this connection we note that infinite-dimensional Lie algebras whose proper subalgebras have dimension $1$ are some analogues of the Tarski Monster in group theory. The problem of the existence of such Lie algebras is one of the most interesting and difficult unsolved problems in the general theory of Lie algebras.

\section{Leibniz algebras, whose subalgebras are either ideals or self-idealizing. The case when the locally nilpotent radical is abelian}

Every Leibniz algebra $L$ possesses the following specific ideal. Denote by $\mathrm{Leib}(L)$ the subspace generated by all elements $[a,a]$, $a\in L$. It is possible to prove that $\mathrm{Leib}(L)$ is an ideal of $L$. The ideal $\mathrm{Leib}(L)$ is called the \textit{Leibniz kernel} of $L$.

We note that the factor-algebra $L/\mathrm{Leib}(L)$ is a Lie algebra, and conversely if $H$ is an ideal of $L$ such that the factor-algebra $L/H$ is a Lie algebra, then $\mathrm{Leib}(L)\leqslant H$.

\begin{lem}\label{L2.1}
Let $L$ be a Leibniz algebra over a field $F$, whose subalgebras are either ideals or self-idealizing. Suppose that $A,B$ are subalgebras of $L$ such that $B$ is an ideal of $A$. If $A/B$ is a non-trivial cyclic algebra, then either $\mathrm{dim}_{F}(A/B)= 1$, or $A/B=Fa_{1} \oplus Fa_{2}$ and $[a_{1},a_{1}]=a_{2}$, $[a_{2},a_{1}]=[a_{2},a_{2}]= 0$, $[a_{1},a_{2}]=\lambda a_{2}$ where $\lambda\in\{0,1\}$.
\end{lem}

\pf Since $B$ is an ideal of $A$, $A\leqslant\mathrm{I}_{L}(B)$, so that $B\neq\mathrm{I}_{L}(B)$, and hence $B$ is an ideal of $L$. Let $a$ be an element of $A$ such that $A/B=\langle a+B\rangle$. If $[a,a]\in B$, then a cyclic subalgebra $\langle a+B\rangle$ has dimension $1$. Suppose now that $d=[a,a]\not\in B$. It follows that $\mathrm{Leib}(A/B)$ is not trivial. This subalgebra is abelian, in particular, every its subspace is a subalgebra. If we suppose that   $\mathrm{dim}_{F}(\mathrm{Leib}(A/B))> 1$, then $\langle d+B\rangle\neq\mathrm{Leib}(A/B)$. In this case, cyclic subalgebra $\langle d+B\rangle$ cannot be self-idealizing. Then $\langle d,B\rangle$ cannot be self-idealizing. It follows that $\langle d,B\rangle$ is an ideal of $L$. Then $\langle d+B\rangle$ is an ideal in factor-algebra $L/B$. In this case, $[a+B,d+B],[d+B,a+B]\in\langle d+B\rangle$. It follows that $\langle a+B\rangle=F(a+B)\oplus F(d+B)$, in particular, $\mathrm{dim}_{F}(\langle a+B\rangle)= 2$. Now we can apply the description of Leibniz algebra of dimension $2$ (see, for example, survey \cite{KKPS2017}).\qed

\begin{lem}\label{L2.2}
Let $L$ be a Leibniz algebra over a field $F$, whose subalgebras are either ideals or self-idealizing. Suppose that $A,B$ are subalgebras of $L$ such that $B$ is an ideal of $A$ and $A/B$ is locally nilpotent. If $A/B$ is not cyclic, then every subalgebra of $A$, including $B$, is an ideal of $L$. In particular, $A$ and $B$ are also ideals of $L$.
\end{lem}

\pf Since $A/B$ is not trivial, then as in the proof of Lemma~\ref{L2.1}, we can prove that $B$ is an ideal of $L$. Let $a$ be an arbitrary element of $A$ such that $a\not\in B$. Since $A/B$ is not cyclic, $\langle a+B\rangle\neq A/B$. Thus we can choose in $A$ an element $b$ such that $b\not\in\langle a,B\rangle$. Since $A/B$ is locally nilpotent, $\langle a+B,b+B\rangle/B=(\langle a,b\rangle+B)/B$ is nilpotent. Then $\mathrm{I}_{(\langle a,b\rangle+B)/B}(\langle a+B\rangle/B)\neq\langle a+B\rangle/B$ \cite{BD2011}. It follows that $\mathrm{I}_{A}(\langle a,B\rangle)\neq\langle a,B\rangle$, so that $\langle a,B\rangle$ is an ideal of $L$.

Let
$$\mathfrak{S}=\{V|\ V\ \mbox{is a subalgebra of $A$ such that $B\leqslant V$ and $V/B$ is a cyclic algebra}\}.$$
If $V\in\mathfrak{S}$ then by proved above, $V$ is an ideal of $L$. Clearly, the subalgebras belonging to $\mathfrak{S}$ generated $A$. It follows that $A$ is an ideal of $L$.\qed

\begin{cor}\label{C2.3}
Let $L$ be a Leibniz algebra over a field $F$, whose subalgebras are either ideals or self-idealizing. Suppose that $A,B$ are subalgebras of $L$ such that $B$ is an ideal of $A$ and $A/B$ is abelian. If $\mathrm{dim}_{F}(A/B)> 1$, then every subalgebra of $A$, including $B$, is an ideal of $L$. In particular, $A$ and $B$ are also ideals of $L$.
\end{cor}

\begin{cor}\label{C2.4}
Let $L$ be a Leibniz algebra over a field $F$, whose subalgebras are either ideals or self-idealizing. If $A$ is an abelian subalgebra of $L$ such that $\mathrm{dim}_{F}(A)> 1$, then every subalgebra of $A$ is an ideal of $L$. In particular, $A$ is also an ideal of $L$.
\end{cor}

\begin{cor}\label{C2.5}
Let $L$ be a Leibniz algebra over a field $F$, whose subalgebras are either ideals or self-idealizing. Suppose that $A$ is a locally nilpotent subalgebra of $L$. If $A$ is not cyclic algebra, then every subalgebra of $A$ is an ideal of $L$. In particular, $A$ is also an ideal of $L$.
\end{cor}

Let $L$ be a Leibniz algebra over a field $F$, $M$ be a non–empty subset of $L$ and $H$ be a subalgebra of $L$. Put
\begin{equation*}
\begin{split}
&\mathrm{Ann}_{H}^{left}(M)=\{a\in H|\ [a,M]=\langle0\rangle\},\\
&\mathrm{Ann}_{H}^{right}(M)=\{a\in H|\ [M,a]=\langle0\rangle\}.
\end{split}
\end{equation*}

The subset  $\mathrm{Ann}_{H}^{left}(M)$ is called the \textit{left annihilator} of $M$ in $H$. The subset $\mathrm{Ann}_{H}^{right}(M)$ is called the \textit{right annihilator} of $M$ in $H$. The intersection
$$\mathrm{Ann}_{H}(M)=\mathrm{Ann}_{H}^{left}(M)\cap\mathrm{Ann}_{H}^{right}(M)=\{a\in H|\ [a,M]=[M,a]=\langle0\rangle\}$$
is called the \textit{annihilator} of $M$ in $H$.

It is not hard to see that all these subsets are the subalgebras of $L$. Moreover, if $M$ is a left ideal of $L$, then it is not hard to prove that $\mathrm{Ann}_{L}^{left}(M)$ is an ideal of $L$. Also it is possible to prove that if $M$ is an ideal of $L$, then $\mathrm{Ann}_{L}^{right}(M)$ is a left ideal of $L$. Furthermore, $\mathrm{Ann}_{L}(M)$ is an ideal of $L$.

The \textit{left} (respectively \textit{right}) \textit{center} $\zeta^{left}(L)$ (respectively $\zeta^{right}(L)$) of a Leibniz algebra $L$ is defined by the rule:
$$\zeta^{left}(L)=\{x\in L|\ [x,y]= 0\ \mbox{for each element }y\in L\}$$
(respectively,
$$\zeta^{right}(L)=\{x\in L|\ [y,x]= 0\ \mbox{for each element }y\in L\}).$$
It is not hard to prove that the left center of $L$ is an ideal of $L$, but it is not true for the right center. Moreover, $\mathrm{Leib}(L)\leqslant\zeta^{left}(L)$, so that $L/\zeta^{left}(L)$ is a Lie algebra. The right center is an subalgebra of $L$, and in general, the left and right centers are different. They even may have different dimensions (see \cite{KOP2016}).

The \textit{center} $\zeta(L)$ of $L$ is defined by the rule:
$$\zeta(L)=\{x\in L|\ [x,y]=[y,x]= 0\ \mbox{for each element }y\in L\}.$$
The center is an ideal of $L$.

Let $L$ be a Leibniz algebra over a field $F$. A linear transformation $f$ of $L$ is called a \textit{derivation}, if
$$f([a,b])=[f(a),b]+[a,f(b)]$$
for all $a,b\in L$. Denote by $\mathrm{End}_{F}(L)$ the set of all linear transformations of $L$, then $L$ is an associative algebra by the operations $+$ and $\circ$. As usual, $\mathrm{End}_{F}(L)$ is a Lie algebra by the operations $+$ and $[-,-]$, where $[f,g]=f\circ g-g\circ f$ for all $f,g\in\mathrm{End}_{F}(L)$. It is possible to show that the subset $\mathrm{Der}(L)$ of all derivations of $L$ is a subalgebra of Lie algebra $\mathrm{End}_{F}(L)$ (see, for example, a survey \cite{KKPS2017}).

\begin{lem}\label{L2.6}
Let $L$ be a Leibniz algebra over a field $F$ and $A$ be an abelian ideal of $L$. If every subalgebra of $A$ is an ideal of $L$, then for each element $x\in L$ there exist the elements $\lambda_{x},\rho_{x}\in F$ such that $[x,a]=\lambda_{x}a$ and $[a,x]=\rho_{x}a$ for every element $a\in A$.
\end{lem}

\pf If $\mathrm{dim}_{F}(A)= 1$, all is trivial. Suppose that $\mathrm{dim}_{F}(A)> 1$. Since $A$ is abelian, a subspace $Fa$ is a subalgebra of $A$ for every element $a\in A$. Thus, a cyclic subalgebra $\langle a\rangle=Fa$ is an ideal of $L$. Let $x\in L$, then $[x,a]=\alpha a$ (respectively $[a,x]=\beta a$) for some elements $\alpha,\beta\in F$. Since $\mathrm{dim}_{F}(A)> 1$, we can choose an element $c\in A$ such that $a$ and $c$ are linearly independent. By the similar reasons, for element $c$, we obtain $[x,c]=\gamma c$ (respectively $[c,x]=\sigma c$) for some elements $\gamma,\sigma\in F$. We have
$$[x,a-c]=[x,a]-[x,c]=\alpha a-\gamma c\ (\mbox{respectively }[a-c,x]=[a,x]-[c,x]=\beta a-\sigma c).$$
On the other hand, $a-c\in A$, so that $F(a-c)$ must be an ideal of $L$, i.e.
$$[x,a-c]=\eta(a-c)=\eta a-\eta c\ (\mbox{respectively }[a-c,x]=\mu(a-c)=\mu a-\mu c)$$
for some $\eta,\mu\in F$. It follows that $\alpha a-\gamma c=\eta a-\eta c$ (respectively $\beta a-\sigma c=\mu a-\mu c$). Hence $\alpha=\eta=\gamma$ (respectively $\beta=\mu=\sigma$).

In other words, for every element $x\in L$ there exist the elements $\lambda_{x},\rho_{x}\in F$ such that $[x,a]=\lambda_{x}a$ and $[a,x]=\rho_{x}a$ for every $a\in A$.\qed

\begin{lem}\label{L2.7}
Let $L$ be a Leibniz algebra over a field $F$, whose subalgebras are either ideals or self-idealizing, and $A$ be a maximal abelian ideal of $L$ including  $\mathrm{Leib}(L)$. Suppose that $\mathrm{dim}_{F}(A)> 1$. Then the following conditions hold:
\begin{enumerate}
\item[\upshape(i)] $A=\zeta^{left}(L)$;
\item[\upshape(ii)] $L/\mathrm{Ann}_{L}(A)$ has dimension $1$;
\item[\upshape(iii)] for each $x\in L$ there exists an element $\sigma_{x}\in F$ such that $[x,a]=\sigma_{x}a$ for every $a\in A$;
\item[\upshape(iv)] every subalgebra of $\mathrm{Ann}_{L}(A)$ is an ideal of $L$.
\end{enumerate}
\end{lem}

\pf For an element $x\in L$ consider the mapping $\mathrm{l}_{x}:A\rightarrow A$, defined by the rule $\mathrm{l}_{x}(a)=[x,a]$ for each element $a\in A$. Then the mapping $\mathrm{l}_{x}$ is a derivation of an ideal $A$, and the set $\{\mathrm{l}_{x}|\ x\in L\}$ is a subalgebra of algebra $\mathrm{Der}(A)$ of all derivations of $A$ (see, for example, a survey \cite{KKPS2017}).

By Corollary~\ref{C2.4}, every subalgebra of $A$ is an ideal of $L$. Then Lemma~\ref{L2.6} implies that for every element $x\in L$ there exist the elements $\sigma_{x},\rho_{x}\in F$ such that $[x,a]=\sigma_{x}a$ and $[a,x]=\rho_{x}a$ for every element $a\in A$.

The inclusion $\mathrm{Leib}(L)\leqslant\zeta^{left}(L)$ implies that $A$ contains an element $a_{0}$ such that $[a_{0},x]=0$ for each element $x\in L$. It follows that $[a,x]=0$ for every element $a\in A$. Since it is true for each element $x\in L$, $A\leqslant\zeta^{left}(L)$. The facts that $\zeta^{left}(L)$ is an abelian ideal of $L$ and $A$ is a maximal abelian ideal of $L$ imply that $A=\zeta^{left}(L)$.

Consider now the mapping $\delta:L\rightarrow F$, defined by the rule $\delta(x)=\sigma_{x}$ for each element $x\in L$. For the elements $x,y\in L$ we have
$$\sigma_{x+y}a=[x+y,a]=[x,a]+[y,a]=\sigma_{x}a+\sigma_{y}a=(\sigma_{x}+\sigma_{y})a,$$
and
$$\sigma_{\beta x}a=[\beta x,a]=\beta[x,a]=\beta(\sigma_{x}a)=(\beta\sigma_{x})a,$$
which shows that $\sigma_{x+y}=\sigma_{x}+\sigma_{y}$ and $\sigma_{\beta x}=\beta\sigma_{x}$ for all $x,y\in L$, $\beta\in F$. It follows that the mapping $\delta$ is linear. Furthermore,
$$\mathrm{Ker}(\delta)=\{x\in L|\ \delta(x)=\sigma_{x}= 0\}.$$
This means that $[x,a]= 0$ for every element $x\in A$. In other words, $\mathrm{Ker}(\delta)\leqslant\mathrm{Ann}_{L}^{left}(A)$. The converse inclusion is obvious, so that $\mathrm{Ker}(\delta)=\mathrm{Ann}_{L}^{left}(A)$. As we remarked above, $\mathrm{Ann}_{L}^{left}(A)$ is a two-side ideal of $L$, so we obtain that $L/\mathrm{Ann}_{L}^{left}(A)$ is isomorphic to $F$, in particular, this factor-algebra has dimension $1$.

The equality $A=\zeta^{left}(L)$ implies that $L=\mathrm{Ann}_{L}^{right}(A)$. This means that $\mathrm{Ann}_{L}^{left}(A)=\mathrm{Ann}_{L}(A)$.

Let $z\in\mathrm{Ann}_{L}(A)$. If $z\in A$, then, as we have noted above, a subalgebra  $\langle z\rangle$ is an ideal of $L$. Assume that $z\not\in A$. By Lemma~\ref{L2.1}, a cyclic subalgebra $\langle a\rangle$ has dimension at most $2$. If $\mathrm{dim}_{F}(\langle z\rangle)= 1$, then $\langle z\rangle=Fz$, and $\langle z\rangle\cap A=\langle0\rangle$. Then for every non-zero element $a\in A$ we have $[a,z]=[z,a]= 0$. This means that $a\in\mathrm{I}_{L}(\langle z\rangle)$, in particular, $\mathrm{I}_{L}(\langle z\rangle)\neq\langle z\rangle$. In this case, $\langle z\rangle$ is an ideal of $L$. Suppose now that $\mathrm{dim}_{F}(\langle z\rangle)= 2$. Put $v=[z,z]$, then $v\in A$. Since $\mathrm{dim}_{F}(A)> 1$, we can choose an element $d\in A$ such that $Fv\cap Fd=\langle0\rangle$. Then $Fd\cap\langle z\rangle=\langle0\rangle$. We have again $[d,z]=[z,d]= 0$. It follows that $d\in\mathrm{I}_{L}(\langle z\rangle)$, so that $\mathrm{I}_{L}(\langle z\rangle)\neq\langle z\rangle$. Hence $\langle z\rangle$ is an ideal of $L$. Thus every cyclic subalgebra of $\mathrm{Ann}_{L}(A)$ is an ideal of $L$. It follows that every subalgebra of $\mathrm{Ann}_{L}(A)$ is an ideal of $L$.\qed

\begin{cor}\label{C2.8}
Let $L$ be a Leibniz algebra over a field $F$, whose subalgebras are either ideals or self-idealizing, and suppose that $\mathrm{dim}_{F}(\mathrm{Leib}(L))> 1$. Then $L/\mathrm{Ln}(L)$ has a dimension at most $1$ and every subalgebra of $\mathrm{Ln}(L)$ is an ideal of $L$.
\end{cor}

\pf Let $A$ be the maximal abelian ideal of $L$ including $\mathrm{Leib}(L)$. Lemma~\ref{L2.7} implies that every subalgebra of $\mathrm{Ann}_{L}(A)$ is an ideal of $L$. Then $\mathrm{Ann}_{L}(A)$ is nilpotent \cite{KSS2017}, so that
$$\mathrm{Ann}_{L}(A)\leqslant\mathrm{Ln}(L).$$
Using Lemma~\ref{L2.7} we obtain that $L/\mathrm{Ln}(L)$ has a dimension at most $1$, and every subalgebra of $\mathrm{Ln}(L)$ is an ideal of $L$.\qed

\medskip

{\noindent\bf Proof of Theorem~A ---} The fact that $\mathrm{Leib}(L)$ is abelian, implies that $\mathrm{Leib}(L)\leqslant A$. Since $A$ is not cyclic, $\mathrm{dim}_{F}(A)> 1$. Being a maximal locally nilpotent ideal, $A$ is a maximal abelian ideal of $L$. By Lemma~\ref{L2.7}, every subalgebra of $A$ is an ideal of $L$, and factor-algebra $L/A$ has dimension at most $1$. Suppose that $L\neq A$. Choose an element $v$ such that $L=A\oplus Fv$. Lemma~\ref{L2.7} shows that $A=\zeta^{left}(L)$. Lemma~\ref{L2.7} also implies that there exists an element $\sigma\in F$ such that $[v,a]=\sigma a$ for each element $a\in A$. Since $A\neq L$, $\sigma\neq 0$. Let $b$ be an arbitrary element of $A$. As we have already noted, the subalgebra $Fb$ is an ideal of $L$. If $d$ is an arbitrary element of $Fb$, then $d=\lambda b$ for some $\lambda\in F$. Then
$$d=\lambda(\sigma^{-1}\sigma)b=\lambda\sigma^{-1}(\sigma b)=\lambda\sigma^{-1}[v,b]=[v,\lambda\sigma^{-1}b]\in[v,Fb].$$
It follows that $[v,Fb]=Fb$. Since it is true for every one dimensional subalgebra of $A$, $A=[v,A]$.

Let $x$ be an arbitrary element of $L$. The fact that $\mathrm{dim}_{F}(L/A)= 1$ implies that $L/A$ is abelian. It follows that $[v,x]=c\in A$. By proved above, $A$ contains an element $u$ such that $c=[v,u]$. Hence $[v,x]=[v,u]$ and therefore $[v,x-u]= 0$. This means that $x-u\in\mathrm{Ann}_{L}^{right}(v)$ or $x\in A+\mathrm{Ann}_{L}^{right}(v)$. Since $x$ is an arbitrary element of $L$, we obtain the equality
$$L=A+\mathrm{Ann}_{L}^{right}(v).$$
Let $a\in A\cap\mathrm{Ann}_{L}^{right}(v)$ and suppose that $a\neq 0$. Then $[v,a]= 0$. On the other hand, $[v,a]=\sigma a$ where $\sigma$ is not zero, and we obtain a contradiction. This contradiction proves that $A\cap\mathrm{Ann}_{L}^{right}(v)=\langle0\rangle$. Put $W=\mathrm{Ann}_{L}^{right}(v)$, then $L=A\oplus W$. Isomorphism $W\cong L/A$ implies that a subalgebra $W$ is abelian and has a dimension $1$, that is, $W=Fw$. We have $w=\lambda v+a_{1}$ for some element $\lambda\in F$ and $a_{1}\in A$. Since $w\not\in A$, $\lambda\neq 0$. Replacing $w$ by $\lambda^{-1}w$, we can suppose that $w=v+a_{2}$. Then $[w,a]=[v,a]=\sigma a$ for each element $a\in A$.

If we suppose that $\mathrm{I}_{L}(W)\neq W$, then a subalgebra $W$ must be an ideal. But in this case, $L$ is abelian, and we obtain a contradiction. This contradiction shows that $W$ is self-idealizing.

Conversely, let $L$ be a Leibniz algebra, satisfying all above conditions. The conditions (i) and (iii) imply that every cyclic subalgebra of $A$ is an ideal of $L$. It follows that every subalgebra of $A$ is an ideal of $L$.

Let $S$ be an arbitrary subalgebra of $L$. If $A$ includes $S$, then by noted above, $S$ is an ideal of $L$. Therefore suppose that $A$ does not include $S$. Then $S$ contains an element $\mu w+e$ where $0 \neq\mu\in F$ and $e\in A$. We have $L=A+S$. If we suppose that $\mathrm{I}_{L}(S)\neq S$, then we can choose an element $a\in A$ such that $[S,a]\leqslant S$ and $a\not\in S$. We have
$$[\mu w+e,a]=\mu[w,a]=\mu\sigma a\in S.$$
Since $\mu\sigma\neq 0$, $a\in S$, and we obtain a contradiction. This contradiction shows that $\mathrm{I}_{L}(S)=S$.\qed

\section{Leibniz algebras, whose subalgebras are either ideals or self-idealizing. The case when the locally nilpotent radical is not abelian.}

\begin{lem}\label{L3.1}
Let $L$ be a Leibniz algebra over a field $F$ and $f$ be a derivation of $L$. Then $f(\zeta^{left}(L))\leqslant\zeta^{left}(L)$, $f(\zeta^{right}(L))\leqslant\zeta^{right}(L)$ and $f(\zeta(L))\leqslant\zeta(L)$.
\end{lem}

\pf Let $x$ be an arbitrary element of $L$ and let $z\in\zeta^{left}(L)$. Then $[z,x]= 0$. Since a derivation is a linear mapping, $f([z,x])= 0$. On the other hand,
$$0 =f(0)=f([z,x])=[f(z),x]+[z,f(x)]=[f(z),x],$$
so that $f(z)\in\zeta^{left}(L)$. Let $z\in\zeta^{right}(L)$. Then $[x,z]= 0$. Now we have
$$0 =f(0)=f([x,z])=[f(x),z]+[x,f(z)]=[x,f(z)],$$
so that $f(z)\in\zeta^{right}(L)$. Both proved above inclusions imply that $f(\zeta(L))\leqslant\zeta(L)$.\qed

\begin{lem}\label{L3.2}
Let $L$ be a Leibniz algebra over a field $F$ where $\mathrm{char}(F)\neq 2$, $L=Fa_{1}\oplus Fa_{2}$ and $[a_{1},a_{1}]=a_{2}$, $[a_{1},a_{2}]=[a_{2},a_{1}]=[a_{2},a_{2}]= 0$. A linear mapping $f$ is a derivation of $L$ if and only if $f(a_{1})=\alpha a_{1}+\beta a_{2}$ for some elements $\alpha,\beta\in F$ and $f(a_{2})= 2\alpha a_{2}$.
\end{lem}

\pf The equality $\zeta(L)=Fa_{2}$ together with Lemma~\ref{L3.1} imply that $f(a_{2})=\gamma a_{2}$ for some element $\gamma\in F$. We have
\begin{equation*}
\begin{split}
\gamma a_{2}&=f(a_{2})=f([a_{1},a_{1}])=[f(a_{1}),a_{1}]+[a_{1},f(a_{1})]\\
&=[\alpha a_{1}+\beta a_{2},a_{1}]+[a_{1},\alpha a_{1}+\beta a_{2}]\\
&=\alpha[a_{1},a_{1}]+\alpha[a_{1},a_{1}]=2\alpha a_{2}.
\end{split}
\end{equation*}
It follows that $f(a_{2})=2\alpha a_{2}$.

Let $f(a_{1})=\alpha a_{1}+\beta a_{2}$, and let $x,y$ be arbitrary elements of $L$. Then $x=\lambda a_{1}+\mu a_{2}$, $y=\sigma a_{1}+\rho a_{2}$. We have
\begin{equation*}
\begin{split}
[x,y]&=[\lambda a_{1}+\mu a_{2},\sigma a_{1}+\rho a_{2}]=\lambda\sigma[a_{1},a_{1}]+\lambda\rho[a_{1},a_{2}]+\mu\sigma[a_{2},a_{1}]+\mu\rho[a_{2},a_{2}]\\
&=\lambda\sigma a_{2};\\
f(x)&=f(\lambda a_{1}+\mu a_{2})=\lambda f(a_{1})+\mu f(a_{2})=\lambda(\alpha a_{1}+\beta a)_{2})+\mu(2\alpha a_{2})\\
&=\lambda\alpha a_{1}+\lambda\beta a_{2}+ 2\mu\alpha a_{2}=\lambda\alpha a_{1}+(\lambda\beta+ 2\mu\alpha)a_{2};\\
f(y)&=f(\sigma a_{1}+\rho a_{2})=\sigma f(a_{1})+\rho f(a_{2})=\sigma(\alpha a_{1}+\beta a_{2})+\rho(2\alpha a_{2})\\
&=\sigma\alpha a_{1}+\sigma\beta a_{2})+ 2\rho\alpha a_{2}=\sigma\alpha a_{1}+(\sigma\beta+ 2\rho\alpha)a_{2};\\
2\alpha\lambda\sigma a_{2}&=\lambda\sigma f(a_{2})=f(\lambda\sigma a_{2})=f([x,y])=[f(x),y]+[x,f(y)]\\
&=[\lambda\alpha a_{1}+(\lambda\beta+ 2\mu\alpha)a_{2},\sigma a_{1}+\rho a_{2}]+[\lambda a_{1}+\mu a_{2},\sigma\alpha a_{1}+(\sigma\beta+ 2\rho\alpha)a_{2}]\\
&=\lambda\alpha\sigma[a_{1},a_{1}]+\lambda\sigma\alpha[a_{1},a_{1}]= 2\alpha\lambda\sigma a_{2}.
\end{split}
\end{equation*}
Thus each linear transformation $f$ of $L$, satisfying $f(a_{1})=\alpha a_{1}+\beta a_{2}$ and $f(a_{2})= 2\alpha a_{2}$ is a derivation of $L$.\qed

Using similar arguments, we obtain

\begin{lem}\label{L3.3}
Let $L$ be a Leibniz algebra over a field $F$ where $\mathrm{char}(F)= 2$, $L=Fa_{1}\oplus Fa_{2}$ and $[a_{1},a_{1}]=a_{2}$, $[a_{1},a_{2}]=[a_{2},a_{1}]=[a_{2},a_{2}]= 0$. A linear mapping $f$ is a derivation of $L$ if and only if $f(a_{2})= 0$.
\end{lem}

\begin{cor}\label{C3.4}
Let $L$ be a Leibniz algebra over a field $F$, where $\mathrm{char}(F)= 2$, $L=Fa_{1}\oplus Fa_{2}$ and $[a_{1},a_{1}]=a_{2}$, $[a_{1},a_{2}]=[a_{2},a_{1}]=[a_{2},a_{2}]= 0$. Then the algebra of derivations of $L$ is isomorphic to a subalgebra of $\mathrm{M}_{2}(F)$, having the following form $\alpha E_{11}+\beta E_{21}$, $\alpha,\beta\in F$. In particular, it is abelian and has dimension $2$.
\end{cor}

\pf If $f$ is a derivation of $L$, then Lemma~\ref{L3.3} shows that $f(a_{1})=\alpha a_{1}+\beta a_{2}$ and $f(a_{2})= 0$ for some elements $\alpha,\beta\in F$. Thus a matrix of a mapping $f$ in a basis $\{a_{1},a_{2}\}$ is $\alpha E_{11}+\beta E_{21}$, $\alpha,\beta\in F$. And conversely, if a linear mapping $f$ has in a basis $\{a_{1},a_{2}\}$ such matrix, then $f$ is a derivation of $L$. It follows that an algebra of derivations of $L$ is isomorphic to subalgebra of matrices, having a form $\alpha E_{11}+\beta E_{21}$, $\alpha,\beta\in F$.\qed

\begin{lem}\label{L3.5}
Let $L$ be a Leibniz algebra over a field $F$, whose subalgebras are either ideals or self-idealizing. Suppose that $A$ is an ideal of $L$ and $B/A$ is a non-trivial subalgebra of $L/A$. If $S/A$ is a subalgebra of $L/A$ such that $S/A\leqslant\mathrm{Ann}_{L/A}(B/A)$ and $S/A$ does not include $B/A$, then $S$ is an ideal of $L$.
\end{lem}

\pf Choose in $B/A$ an element $zA$ such that $zA\not\in S/A$. The inclusion $S/A\leqslant\mathrm{Ann}_{L/A}(B/A)$ implies that $[zA,S/A]=[S/A,zA]=\langle0\rangle$. It follows that $[z,S],[S,z]\leqslant A\leqslant S$. The choice of an element $z$ implies that $z\not\in S$. This means that $\mathrm{I}_{L}(S)\neq S$. It follows that $S$ is an ideal of $L$.\qed

\begin{lem}\label{L3.6}
Let $L$ be a Leibniz algebra over a field $F$, whose subalgebras are either ideals or self-idealizing. Suppose that $L\neq K=\mathrm{Ln}(L)=K$ is not cyclic and not abelian. Then every subalgebra of $K$ is an ideal of $L$. If $[K,K]\leqslant\zeta(L)$, then $K$ is a strong extraspecial algebra, $[K,K]=\zeta(K)=\mathrm{Leib}(L)\leqslant\zeta(L)$.
\end{lem}

\pf Using Corollary~\ref{C2.5} we obtain that every subalgebra of $K$ is an ideal of $L$. Since $K$ is non-abelian, by above noted $K=E\oplus Z$ where $Z$ is a subalgebra of the center of $K$ and $E$ is a strong extraspecial algebra. Since $K$ is non-abelian, $E$ is not trivial. Choose in $E$ an element $y$ such that $z=[y,y]\neq 0$. Let $Y$ be a subalgebra, generated by element $y$, then $Y=Fy\oplus Fz$ and $[z,y]=[y,z]= 0$. By above remarked a subalgebra $Y$ is an ideal of $L$. Since $[K,K]=[E,E]$ has dimension $1$, the fact that $z\in[E,E]$ implies that $Fz=[K,K]$. It follows that $\in\zeta(L)$.

Suppose that $Z$ is not trivial and consider a subalgebra $\langle z\rangle\oplus Z$. By imposed conditions, this subalgebra and every its subalgebra are ideals of $L$. Then Lemma~\ref{L2.6} implies that $Z\leqslant\zeta(L)$. Consider now $L/\langle z\rangle$. Its ideal $K/\langle z\rangle$ is abelian and, clearly, $\mathrm{dim}_{F}(K/\langle z\rangle)> 1$. Lemma~\ref{L2.2} implies that every subalgebra of $K/\langle z\rangle$ is an ideal of $L/\langle z\rangle$. Using the fact that the factor $(\langle z\rangle\oplus Z)/\langle z\rangle$ is central in $L$ and Lemma~\ref{L2.6} we obtain that the factor $K/\langle z\rangle$ is central in $L/\langle z\rangle$. Since $L\neq K$, we can choose an element $v$ such that $v\not\in K$. By Lemma~\ref{L2.1} a subalgebra $V=\langle v\rangle$ has dimension $1$ or $2$. Suppose that $\mathrm{dim}_{F}(V)= 1$, then $V\cap K=\langle0\rangle$, in particular, $z\not\in V$ but $z\in\mathrm{I}_{L}(V)$. It follows that $V$ is an ideal of $L$. The fact that $\mathrm{dim}_{F}(V)= 1$ implies that $V$ is abelian. But in this case $\mathrm{Ln}(L)=K$ must includes $V$, and we obtain a contradiction. Suppose now that $\mathrm{dim}_{F}(V)= 2$. If we assume that $[V,V]\cap K=\langle0\rangle$, then again $z\not\in[V,V]$ and $z\in\mathrm{I}_{L}([V,V])$. It follows that $[V,V]$ is an ideal of $L$. The fact that $\mathrm{dim}_{F}([V,V])= 1$ implies that $[V,V]$ is abelian. But in this case $\mathrm{Ln}(L)=K$ must includes $[V,V]$, and we obtain a contradiction with $[V,V]\cap K=\langle0\rangle$. Hence $[V,V]\cap K\neq\langle0\rangle$. If $\langle0\rangle\neq[V,V]\cap\langle z\rangle$, then the fact that $\mathrm{dim}_{F}([V,V])= 1$ shows that $[V,V]=\langle z\rangle$. Then $\langle0\rangle=V/\langle z\rangle\cap K/\langle z\rangle$. Since $K/\langle z\rangle$ is central in $L/\langle z\rangle$, $\mathrm{I}_{L}(V)\neq V$, so that $V$ is an ideal of $L$. The fact that $[V,V]=\langle z\rangle\leqslant\zeta(L)$ implies that $V$ is nilpotent. But in this case again $K$ must includes $V$, and we obtain a contradiction. This contradiction shows that $\langle 0\rangle=[V,V]\cap\langle z\rangle$. Then $V\cap\langle z\rangle=\langle0\rangle$. In this case $z\not\in V$ but $z\in\mathrm{I}_{L}(V)$. It follows that $V$ is an ideal of $L$. On the other hand
$$([V,V]+\langle z\rangle)/\langle z\rangle\leqslant K/\langle z\rangle\leqslant\zeta(L/\langle z\rangle).$$
It follows that $(V+\langle z\rangle)/\langle z\rangle$ is nilpotent. An inclusion $\langle z\rangle\leqslant\zeta(L)$ implies that $V+\langle z\rangle$ is nilpotent. In particular, $V$ itself is nilpotent. Then again $K$ must includes $V$, and we obtain a final contradiction, which proves that $Z=\langle0\rangle$, i.e. $K=E$ is a strong extraspecial Leibniz algebra.

Suppose now that $\mathrm{Leib}(L)\neq\langle z\rangle$. Since $\mathrm{Leib}(L)$ is an abelian ideal, $K$ includes $\mathrm{Leib}(L)$. Then $\mathrm{dim}_{F}(\mathrm{Leib}(L))> 1$, so that $\mathrm{Leib}(L)$ must contains an element $u\not\in\langle z\rangle$. Since $\mathrm{Leib}(L)$ is abelian, $[u,u]= 0$. On the other hand, $K$ is a strong extraspecial Leibniz algebra and it follows that $[u,u]\neq 0$. This contradiction proves that $\mathrm{Leib}(L)=\langle z\rangle$.\qed

\begin{cor}\label{C3.7}
Let $L$ be a Leibniz algebra over a field $F$, whose subalgebras are either ideals or self-idealizing. Suppose that $L\neq\mathrm{Ln}(L)=K$ is not cyclic and not abelian. If $[K,K]\leqslant\zeta(L)$, then $\mathrm{char}(F)= 2$.
\end{cor}

\pf By Lemma~\ref{L3.6}, $K$ is a strong extraspecial algebra, $[K,K]=\zeta(K)=\mathrm{Leib}(L)$. Since $K$ is non-abelian, $E$ is not trivial. Choose in $E$ an element $y$ such that $z=[y,y]\neq 0$. Let $Y$ be a subalgebra, generated by $y$, then $Y=Fy\oplus Fz$ and $[z,y]=[y,z]= 0$. By our conditions, $Y$ is an ideal of $L$. As in Lemma~\ref{L3.6} we can show that $\langle z\rangle=[K,K]$, so that $z\in\zeta(L)$.

For an element $x\in L$ consider the mapping $\mathrm{l}_{x}:Y\rightarrow Y$, defined by the rule $\mathrm{l}_{x}(a)=[x,a]$ for each element $a\in Y$. Then the mapping $\mathrm{l}_{x}$ is a derivation of $Y$.

Suppose now that $\mathrm{char}(F)\neq 2$. Using Lemma~\ref{L3.2} and the fact that $\mathrm{l}_{x}(z)= 0$ we obtain that $\mathrm{l}_{x}(y)=\beta z$ for some element $\beta\in F$. In other words, $[x,y]\in\langle z\rangle$ for all elements $x\in L$. It follows that
$$Y/\langle z\rangle\leqslant\zeta^{right}(L/\langle z\rangle).$$
By Lemma~\ref{L3.6} $\mathrm{Leib}(L)=\langle z\rangle$. It follows that $L/\langle z\rangle$ is a Lie algebra. Then we have $\zeta^{right}(L/\langle z\rangle)=\zeta(L/\langle z\rangle)$, so that $Y/\langle z\rangle\leqslant\zeta(L/\langle z\rangle)$. Lemma~\ref{L2.2} implies that every subalgebra of $K/\langle z\rangle$ is an ideal of $L/\langle z\rangle$. Using Lemma~\ref{L2.6}, we obtain that the factor $K/\langle z\rangle$ is central in $L/\langle z\rangle$. Since $L/\langle z\rangle$ is a Lie algebra, every cyclic subalgebra $X/\langle z\rangle$ of $L/\langle z\rangle$ has dimension $1$. It follows that either $K/\langle z\rangle$ includes $X/\langle z\rangle$ or $K/\langle z\rangle\cap X/\langle z\rangle=\langle0\rangle$. In the first case, Lemma~\ref{L3.6} shows that $X$ is an ideal of $L$. In the second case, Lemma~\ref{L3.5} shows that $X$ is an ideal of $L$. Thus, every cyclic subalgebra of $L/\langle z\rangle$ is an ideal of $L/\langle z\rangle$. Since $L/\langle z\rangle$ is a Lie algebra, $L/\langle z\rangle$ is abelian. It follows that $L$ is nilpotent, and we obtain a contradiction with $L\neq K$. This contradiction shows that $\mathrm{char}(F)= 2$.\qed

\medskip

{\noindent\bf Proof of Theorem~B1 ---} Corollary~\ref{C3.7} implies that $\mathrm{char}(F)= 2$. Denote by $K$ the locally nilpotent radical of $L$. By Lemma~\ref{L3.6}, every subalgebra of $K$ is an ideal of $L$ and $K$ is a strong extraspecial algebra. Choose in $K$ an element $y$ such that $z=[y,y]\neq 0$. Let $Y$ be a subalgebra, generated by element $y$, then $Y=Fy\oplus Fz$ and $[z,y]=[y,z]= 0$. By above noted, $Y$ is an ideal of $L$. Lemma~\ref{L3.6} implies also that $Z=Fz=\langle z\rangle=[K,K]=\mathrm{Leib}(L)\leqslant\zeta(L)$. It follows that factor–algebra $L/Z$ is non-abelian (otherwise $L$ must be nilpotent, but this contradicts $\mathrm{Ln}(L)\neq L$). Equality $Z=\mathrm{Leib}(L)$ implies that $L/Z$ is a Lie algebra.

Suppose that $L$ contains an element $w\not\in Z$ such that $[w,w]= 0$. Then $\langle w\rangle=Fw$ and $\langle w\rangle\cap Z=\langle0\rangle$. Inclusion $Z\leqslant\zeta(L)$ implies that $z\in\mathrm{I}_{L}(\langle w\rangle)$, in particular, $\mathrm{I}_{L}(\langle w\rangle)\neq\langle w\rangle$. It follows that $\langle w\rangle$ is an ideal of $L$. Since $\langle w\rangle$ is abelian, $\langle w\rangle\leqslant\mathrm{Ln}(L)$. On the other hand, as we have noted above, $L$ is a strong extraspecial algebra, and we obtain a contradiction. This contradiction shows that $L$ satisfies condition (i).

Since $K$ is not cyclic, $K\neq Y$. It follows that $\mathrm{dim}_{F}(K/Z)> 1$. Let $A=\mathrm{Ann}_{L}^{left}(Y)$. Since $Y$ is an ideal of $L$, $A$ also is an ideal of $L$. Inclusion $Z\leqslant\zeta(L)$ implies that $Z\leqslant A$. Clearly $y\not\in A$, so that $A\cap Y=Z$. It follows that the intersection $Y/Z\cap A/Z$ is trivial. Let $a$ be an arbitrary element of $A$ such that $a\not\in Z$. Then $a\not\in Y$ and $[a,y]=[y,a]= 0$. It follows that $\mathrm{I}_{L}(\langle a\rangle)\neq\langle a\rangle$ and therefore $\langle a\rangle$ is an ideal of $L$. Since $[a,a]\in Z$, $\langle a\rangle$ is nilpotent, so that $\langle a\rangle\leqslant\mathrm{Ln}(L)$. Thus $A\leqslant\mathrm{Ln}(L)$. By Proposition~3.2 of \cite{KOP2016} and Corollary~\ref{C3.4} $\mathrm{dim}_{F}(L/A)\leqslant 2$. Then
$$\mathrm{dim}_{F}(L/(A+Y))\leqslant 1.$$
Since $(A+Y)\leqslant\mathrm{Ln}(L)$ and $L\neq\mathrm{Ln}(L)$, we obtain that
$$\mathrm{dim}_{F}(L/\mathrm{Ln}(L))= 1.$$
Choose an element $u$ such that $u\not\in K$. Then $L=K\oplus Fu$. Corollary~\ref{C2.3} implies that every subalgebra of $K$, including $Z$, is an ideal of $L$. Lemma~\ref{L2.6} implies that there exists an element $\alpha\in F$ such that $[a+Z,u+Z]=\alpha(a+Z)$ for each element $a\in K$. We note that $\alpha\neq 0$. Then $[a,u]=\alpha a+z_{a}$ for some element $z_{a}\in Z$. Let $a,b$ be elements of $K$ such that $a,b\not\in Z$. Since $[a,b]\in Z$, $[[a,b],v]= 0$. We have
\begin{equation*}
\begin{split}
0&=[[a,b],u]=[a,[b,u]]-[b,[a,u]]=[a,\alpha b+z_{b}]-[b,\alpha a+z_{a}]\\
&=\alpha[a,b]-\alpha[b,a]=\alpha([a,b]-[b,a]).
\end{split}
\end{equation*}
Since $\alpha\neq 0$, $[a,b]=[b,a]$.

By Lemma~\ref{L3.2} $[u,y]=\lambda y+\mu z$. Using now the fact that every subalgebra of $K/Z$ is an ideal of $L/Z$ and Lemma~\ref{L2.6} we obtain that $[u+Z,a+Z]=\lambda a+Z$ for each element $a\in K\setminus Z$. Since $u\not\in K$, $\lambda\neq 0$. Then put $v=\lambda^{-1}u$ and we will have $[v+Z,a+Z]=a+Z$ for each element $a\in K\setminus Z$.

Consider now a cyclic subalgebra $\langle v\rangle$. By Lemma~\ref{L2.1}, either $\langle v\rangle=Fv$ or $\langle v\rangle=Fv\oplus F[v,v]$. In the first case, $\langle v\rangle\cap\mathrm{Leib}(L)=\langle0\rangle$. Inclusion $\mathrm{Leib}(L)\leqslant\zeta(L)$ shows that $z\in\mathrm{I}_{L}(\langle v\rangle)$. Then $\langle v\rangle$ must be an ideal of $L$. In this situation we have $\langle v\rangle\cap K=\langle0\rangle$. But then $v\in\mathrm{Ann}_{L}(K)$, and we obtain a contradiction. This contradiction shows that $[v,v]=\eta z$ for some non-zero element $\eta\in F$.\qed

We note some details of the structure of $\mathrm{Ln}(L)$.

Let $\{y+Z,w_{\lambda}+Z|\ \lambda\in\Lambda\}$ be a basis of $K/Z$. Put $a_{1}=y$. Since $K/Z$ is abelian, $[w_{\lambda},a_{1}]=\xi_{\lambda}z$ for some element $\xi_{\lambda}\in F$, $\lambda\in\Lambda$. If $\xi_{\lambda}= 0$, then put $u_{\lambda}=w_{\lambda}$. If $\xi_{\lambda}\neq 0$, then put $u_{\lambda}=\xi\lambda a_{1}-w_{\lambda}$. Then
$$[u_{\lambda},a_{1}]=[\xi_{\lambda}a_{1}-w_{\lambda},a_{1}]=\xi_{\lambda}[a_{1},a_{1}]-[w_{\lambda},a_{1}]=\xi_{\lambda}z-\xi_{\lambda}z= 0.$$
Equality $[u_{\lambda},a_{1}]=[a_{1},u_{\lambda}]$ implies that $[a_{1},u_{\lambda}]= 0$. Clearly, elements $\{a_{1}+Z,u_{\lambda}+Z|\ \lambda\in\Lambda\}$ form the basis of $K/Z$. Let $U_{1}/Z$ be the subspace of $K/Z$, generated by elements $\{u_{\lambda}+Z|\ \lambda\in\Lambda\}$. The fact that $K/Z$ is abelian implies that $U_{1}$ is a subalgebra of $K$, and $[U_{1},a_{1}]=[a_{1},U_{1}]=\langle0\rangle$.

Using the similar arguments and transfinite induction, we construct the basis $\{a_{\mu}+Z|\ \mu\in\mathrm{M}\}$ such that $[a_{\mu},a_{\nu}]=[a_{\nu},a_{\mu}]= 0$ for all $\mu,\nu\in\mathrm{M}$, $\mu\neq\nu$. Furthermore $[v,a_{\mu}]=a_{\mu}+\gamma_{\mu}z$ for all $\mu\in\mathrm{M}$.

\medskip

{\noindent\bf Proof of Theorem~B2 ---} Using Corollary~\ref{C2.5} we obtain that every subalgebra of $K$ is an ideal of $L$. Since $K$ is non-abelian, by above noted $K=E\oplus Z$ where $Z$ is a subalgebra of the center of $K$ and $E$ is a strong extraspecial algebra. Since $K$ is non-abelian,$E$ is not trivial. Choose in $E$ an element $y$ such that $z=[y,y]\neq 0$. Let $Y$ be a subalgebra, generated by $y$, then $Y=Fy\oplus Fz$ and $[z,y]=[y,z]= 0$. By our conditions, a subalgebra $Y$ is an ideal of $L$. Since $[K,K]=[E,E]$ has dimension $1$, the fact that $z\in[E,E]$ implies that $Fz=[K,K]$. It follows also that $z\in\zeta(K)$.

For the element $x\in L$ consider the mapping $\mathrm{l}_{x}:Y\rightarrow Y$, defined by the rule $\mathrm{l}_{x}(a)=[x,a]$ for each element $a\in Y$. Then the mapping $\mathrm{l}_{x}$ is a derivation of an ideal $Y$. Since $z\not\in\zeta(L)$, Lemmas~\ref{L3.2} and \ref{L3.3} shows that $\mathrm{char}(F)\neq 2$. Moreover, Lemma~\ref{L3.2} shows that there exist an element $u\in A$ such that $[u,y]=\alpha y+\beta z$ and $[u,z]= 2\alpha z$ and $\alpha\neq 0$. Put $v=\alpha^{-1}u$, then $[v,y]=y+\gamma z$, $[v,z]= 2z$ where $\gamma=\alpha^{-1}\beta$.

Suppose that $Z$ is not trivial and consider $\langle z\rangle\oplus Z$. By imposed conditions, this subalgebra and every its subalgebra are ideals of $L$. Then Lemma~\ref{L2.6} implies that $[v,w]= 2w$ for each element $w\in Z$. Consider now factor-algebra $L/\langle z\rangle$. Its ideal $K/\langle z\rangle$ is abelian and, clearly, $\mathrm{dim}_{F}(K/\langle z\rangle)> 1$. Since $Z\neq\langle z\rangle$, there exists an element $w\in Z$ such that $\langle w\rangle\cap Y=\langle0\rangle$. Lemma~\ref{L2.2} implies that every subalgebra of $K/\langle z\rangle$ is an ideal of $L/\langle z\rangle$. We have
$$[v+\langle z\rangle,y+\langle z\rangle]=[v,y]+\langle z\rangle=y+\langle z\rangle.$$
Using now Lemma~\ref{L2.6} we obtain that $[v+\langle z\rangle,w+\langle z\rangle]=w+\langle z\rangle$. On the other hand,
$$[v+\langle z\rangle,w+\langle z\rangle]=[v,w]+\langle z\rangle= 2w+\langle z\rangle,$$
and we obtain a contradiction. This contradiction proves an equality $Z=\langle0\rangle$, i.e. $K=E$ is a strong extraspecial subalgebra.

Suppose now that $\mathrm{Leib}(L)\neq\langle z\rangle$. Since $\mathrm{Leib}(L)$ is an abelian ideal, $K$ includes $\mathrm{Leib}(L)$. Then $\mathrm{dim}_{F}(\mathrm{Leib}(L))> 1$, so that $\mathrm{Leib}(L)$ must contains an element $u\not\in\langle z\rangle$. Since $\mathrm{Leib}(L)$ is abelian, $[u,u]= 0$. On the other hand, $K$ is a strong extraspecial Leibniz algebra and it follows that $[u,u]\neq 0$. This contradiction proves that $\mathrm{Leib}(L)=\langle z\rangle$.

Since $\langle z\rangle$ is an ideal, having dimension $1$, $L/\mathrm{Ann}_{L}^{left}(\langle z\rangle)\cong F$, in particular, $\mathrm{dim}_{F}(L/\mathrm{Ann}_{L}^{left}(\langle z\rangle)= 1$. The fact that $v\not\in\mathrm{Ann}_{L}^{left}(\langle z\rangle)$ implies that $L=\mathrm{Ann}_{L}^{left}(\langle z\rangle)+\langle v\rangle$. Since the left center of $L$ includes the Leibniz kernel, $L=\mathrm{Ann}_{L}^{right}(\langle z\rangle)$, so that
$$\mathrm{Ann}_{L}(\langle z\rangle)=\mathrm{Ann}_{L}^{left}(\langle z\rangle)\cap\mathrm{Ann}_{L}^{right}(\langle z\rangle)=\mathrm{Ann}_{L}^{left}(\langle z\rangle)\cap L=\mathrm{Ann}_{L}^{left}(\langle z\rangle).$$

Let $x\in\mathrm{Ann}_{L}(\langle z\rangle)$, then Lemma~\ref{L3.2} shows that $[x,y]\in\langle z\rangle$. Then Lemma~\ref{L3.5} implies that $\langle x\rangle$ is an ideal of $L$. An equality $\mathrm{Leib}(L)=\langle z\rangle$ implies that $L/\langle z\rangle$ is a Lie algebra. It follows that $[x+\langle z\rangle,x+\langle z\rangle]=\langle z\rangle$, so that either $\langle x\rangle=Fx$ or $\langle x\rangle=Fx+Fz$. In the second case $[x,x]\in\langle z\rangle$. An equality $[x,z]= 0$ implies that $\langle x\rangle$ is nilpotent. Since this subalgebra is an ideal, $\mathrm{Ln}(L)=K$ includes $\langle x\rangle$. Hence $\mathrm{Ann}_{L}(\langle z\rangle)\leqslant K$ and the fact that $z\in\zeta(K)$ implies that $K=\mathrm{Ann}_{L}(\langle z\rangle)$. It follows that $L=K+\langle v\rangle$.

We have proved above that $[v,x]=x+\xi_{x}z$ for each $x\in K\setminus\langle z\rangle$ and some $\xi_{x}\in F$. Since $L/\langle z\rangle$ is a Lie algebra, then $[v,v]\in\langle z\rangle$, so that $[v,v]=\nu z$ for some $\nu\in F$ and $[v,z]= 2z$.\qed

\begin{ex}\label{E3.8}
Let $L$ be a Leibniz algebra over a field $F$, $\mathrm{char}(F)= 2$, $L$ generated by the elements $a,b,v$ such that
\begin{equation*}
\begin{split}
&[a,a]=z, [b,b]=\sigma z, [v,v]=\eta z,\\
&[z,z]=[z,a]=[a,z]=[z,b]=[b,z]=[z,d]=[d,z]= 0,\\
&[a,b]=[b,a]= 0, [a,v]=a, [v,a]=a+z,\\
&[b,v]=b, [v,b]=b+z.
\end{split}
\end{equation*}
Moreover, the polynomials $X^{2}+\sigma$ and $X^{2}+\eta$ have no root in $F$. It is possible to check that $L$ is a Leibniz algebra, whose subalgebras are either ideals or self-idealizing.
\end{ex}

We say that a field $F$ is \textit{$2$-closed} if a polynomial $X^{2}+\alpha\in F[X]$ has a root in a field $F$ for every non-zero element $\alpha\in F$. This means that a multiplicative group $\mathrm{U}(F)$ of a field $F$ is $2$-divisible.

In particular, every finite field $F$ of characteristic $2$ is $2$-closed. Indeed, $|F|= 2^{n}$ for some positive integer $n$, so that $|\mathrm{U}(F)|=2^{n}- 1$ is odd. It follows that $\mathrm{U}(F)$ is $2$-divisible. As a corollary we obtain that every locally finite field $F$ of characteristic $2$ is $2$-closed.

\begin{cor}\label{C3.9}
Let $L$ be a Leibniz algebra over a field $F$, whose subalgebras are either ideals or self-idealizing. Suppose that $\mathrm{char}(F)= 2$, $L\neq\mathrm{Ln}(L)$ and $\mathrm{Ln}(L)$ is not abelian. If $F$ is $2$-closed, then $\mathrm{Ln}(L)$ is cyclic.
\end{cor}

\pf Denote by $K$ the locally nilpotent radical of $L$ and suppose that $K$ is not cyclic. By Lemma~\ref{L3.6} every subalgebra of $K$ is an ideal of $L$ and $K$ is a strong extraspecial algebra. Choose in $K$ an element $y$ such that $z=[y,y]\neq 0$. Let $Y$ be a subalgebra, generated by element $y$, then $Y=Fy\oplus Fz$ and $[z,y]=[y,z]= 0$. By above noted, $Y$ is an ideal of $L$. Lemma~\ref{L3.6} also implies that $Z=Fz=\langle z\rangle=\mathrm{Leib}(L)\leqslant\zeta(L)$. Since $K$ is not cyclic, $K\neq Y$. Using the argument from the proof of Theorem~B we can find an element $a\in K$ such that $a\not\in Z$, $[a,y]=[y,a]= 0$. We have $[a,a]=\gamma z$ for some element $\gamma\in F$. Let $b=\lambda a+\mu y$, $\lambda,\mu\in F$. We have
\begin{equation*}
\begin{split}
[b,b]&=[\lambda a+\mu y,\lambda a+\mu y]=\lambda^{2}[a,a]+\lambda\mu[a,y]+\mu\lambda[y,a]+\mu^{2}[y,y]\\
&=\lambda^{2}\gamma z+\mu^{2}z=(\lambda^{2}\gamma+\mu^{2})z.
\end{split}
\end{equation*}
Since $F$ is $2$-closed, there exists an element $\sigma\in F$ such that $\sigma^{2}=\gamma$. Put $\lambda= 1$, $\mu=\sigma$, then $a+\sigma y\not\in Z$, $[a+\sigma y,a+\sigma y]=(\gamma+\sigma^{2})z=(\gamma+\gamma)z= 0$. Thus we obtain a contradiction with the fact, that $K$ is a strong extraspecial algebra. This contradiction shows that $K$ must be cyclic.\qed

As we can see now, the next natural step is the consideration of the case when the locally nilpotent radical is cyclic. Lemma~\ref{L2.1} shows that in this case it has a dimension $1$ or $2$.

\medskip

{\noindent\bf Proof of Theorem~C ---} Denote by $K$ the locally nilpotent radical of $L$. By Lemma~\ref{L2.1}, $K$ has a basis $\{a,z\}$ such that $[a,a]=z$, $[z,a]=[a,z]= 0.$

For the element $x\in L$ consider the mapping $\mathrm{l}_{x}:K\rightarrow K$, defined by the rule $\mathrm{l}_{x}(y)=[x,y]$ for each element $y\in K$. Then the mapping $\mathrm{l}_{x}$ is a derivation of $K$.  By Lemmas~\ref{L3.2} and \ref{L3.3}, $\mathrm{l}_{x}(z)= 0$. In other words, $[x,z]= 0$ for all elements $x\in L$. It follows that $z\in\zeta^{right}(L)$. On the other hand, $z\in\mathrm{Leib}(L)$, and inclusion $\mathrm{Leib}(L)\leqslant\zeta^{left}(L)$ implies that $z\in\zeta^{left}(L)$. Hence $z\in\zeta(L)$. We note that $\mathrm{Leib}(L)$ is abelian, therefore $\mathrm{Leib}(L)\leqslant K$. Since $Z=\langle z\rangle$ is a maximal abelian subalgebra of $K$, $Z=\mathrm{Leib}(L)$.

Let $A=\mathrm{Ann}_{L}^{left}(K)$. Since $K$ is an ideal of $L$, $A$ also is an ideal of $L$. Equality $Z=\zeta(K)$ implies that $Z\leqslant A$. Clearly, $a\not\in A$, so that $A\cap K=Z$. Suppose that $K$ does not include $A$, then $A/Z$ is not trivial and the intersection $K/Z\cap A/Z$ is trivial. Let $d$ be an arbitrary element of $A$ such that $d\not\in Z$. Then $d\not\in K$ and $[a,d]=[d,a]= 0$. It follows that $\mathrm{I}_{L}(\langle d\rangle)\neq\langle d\rangle$ and, therefore, $\langle d\rangle$ must be an ideal of $L$. Since $[d,d]\in\mathrm{Leib}(L)=Z$, $\langle d\rangle$ is nilpotent, so that $\langle d\rangle\leqslant K$, and we obtain a contradiction. This contradiction proves the inclusion $A\leqslant K$.

Suppose that $\mathrm{char}(F)\neq 2$. Lemma~\ref{L3.2} shows that in this case, algebra of derivation of $K$ has dimension $1$. Using Proposition~3.2 of \cite{KOP2016} we obtain that $\mathrm{dim}_{F}(L/A)\leqslant 1$. Equality $A=Z$ imply that $L=K$.

Suppose now that $\mathrm{char}(F)= 2$. By Proposition~3.2 of \cite{KOP2016}, Corollary~\ref{C3.4} and equality $A=Z$ imply that $\mathrm{dim}_{F}(L/Z)\leqslant 2$. If $\mathrm{dim}_{F}(L/Z)= 1$, then $L=K$. If $\mathrm{dim}_{F}(L/Z)= 2$, then $L/Z$ is a Lie algebra of dimension $2$. Since $L\neq K$, we can choose an element $v\not\in K$ such that $[a+Z,v+Z]=[v+Z,a+Z]=a+Z$. It follows that $[v,a]=a+\lambda z$, $[a,v]=a+\mu z$ for some elements $\lambda,\mu\in F$.

Finally, consider cyclic subalgebra $\langle v\rangle$. By Lemma~\ref{L2.1}, either $\langle v\rangle=Fv$ or $\langle v\rangle=Fv\oplus F[v,v]$. In the first case, $\langle v\rangle\cap\mathrm{Leib}(L)=\langle0\rangle$. Equality $\mathrm{Leib}(L)=\zeta(L)$ shows that $z\in\mathrm{I}_{L}(\langle v\rangle)$. Then $\langle v\rangle$ must be an ideal of $L$. In this situation we have $\langle v\rangle\cap K=\langle0\rangle$. But then $v\in\mathrm{Ann}_{L}(K)$, and we obtain a contradiction. This contradiction shows that $[v,v]=\eta z$ for some non-zero element $\eta\in F$.

Let $b=\sigma a+\tau v$, $\sigma,\tau\in F$. We have
\begin{equation*}
\begin{split}
[b,b]&=[\sigma a+\tau v,\sigma a+\tau v]=\sigma^{2}[a,a]+\sigma\tau[a,v]+\tau\sigma[v,a]+\tau^{2}[v,v]\\
&=\sigma^{2}z+\sigma\tau(a+\mu z)+\tau\sigma(a+\lambda z)+\tau^{2}\eta z\\
&=\sigma^{2}z+\sigma\tau a+\sigma\tau\mu z+\tau\sigma a+\tau\sigma\lambda z+\tau^{2}\eta z\\
&=(\sigma^{2}+\sigma\tau\mu+\tau\sigma\lambda+\tau^{2}\eta)z=\tau^{2}(\sigma^{2}\tau^{-2}+\sigma\tau^{-1}(\mu+\lambda)+\eta)z.
\end{split}
\end{equation*}
If we suppose that a polynomial $X^{2}+(\mu+\lambda)X+\eta$ has a root $\gamma$ in $F$, then putting $\tau= 1$, $\sigma=\gamma$, we obtain that $[b,b]= 0$, and therefore $\langle b\rangle=Fb$. Clearly, $b\not\in K$. As above $z\in\mathrm{I}_{L}(\langle b\rangle)$, so that $\langle b\rangle$ must be an ideal. Being abelian, $\langle b\rangle$ lies in the locally nilpotent radical of $L$, and we obtain a contradiction. This contradiction shows that polynomial $X^{2}+(\mu+\lambda)X+\eta$ has no root in $F$.\qed

We can show an example of Leibniz algebra satisfying the conditions of Theorem~C.

\begin{ex}\label{E3.10}
Let $L$ be a Leibniz algebra over a field $F$, where $\mathrm{char}(F)= 2$ and $F$ is not $2$-closed. Choose in $F$ an element $\eta$ such that a polynomial $X^{2}+\eta$ have no root in $F$. Let $L$ be a vector space over $F$ and $\{z,a,v\}$ be a basis of $L$. Define the operation $[-,-]$ in the following way:
$$[z,z]=[z,a]=[a,z]=[z,v]=[v,z]= 0,[a,a]=z,[v,v]=\eta z,[v,a]=[a,v]=a.$$
It is possible to check that $L$ is a Leibniz algebra and every subalgebra of $L$ which is not an ideal is self-idealizing.
\end{ex}

\begin{prop}\label{P3.11}
Let $L$ be a Leibniz algebra over a field $F$, whose subalgebras are either ideals or self-idealizing. If $\mathrm{Ln}(L)$ has dimension $1$, then the factor-algebra $L/\mathrm{Ln}(L)$ does not include abelian subalgebras of dimension $2$.
\end{prop}

\pf Denote by $K$ the locally nilpotent radical of $L$. Since $\mathrm{Leib}(L)$ is an abelian ideal, $\mathrm{Leib}(L)\leqslant K$. It follows that $\mathrm{Leib}(L)=K$, so that $L/K$ is a Lie algebra. Let $A=\mathrm{Ann}_{L}^{left}(K)$. Since $K$ is an ideal of $L$, $A$ also is an ideal of $L$. The fact that $\mathrm{dim}_{F}(K)= 1$ implies that $\mathrm{dim}_{F}(L/A)= 1$.

Suppose the contrary, let $U/K$ be the abelian subalgebra of $L/K$ such that $\mathrm{dim}_{F}(U/K)> 1$. Then Corollary~\ref{C2.3} implies that every subalgebra of $U$ including $K$, is an ideal of $L$. Equality $\mathrm{dim}_{F}(L/A)= 1$ shows that $U\cap A> K$. Since $(U\cap A)/K$ is abelian, an ideal $U\cap A$ is nilpotent. Then $U\cap A\leqslant K$, and we obtain a contradiction. This contradiction proves that every abelian subalgebra of $L/K$ has dimension 1.\qed

Thus we can see that the study of Leibniz algebras, whose subalgebras which are not ideals are self-idealizing, reduces to the study of Lie algebras, whose abelian subalgebras of dimension $1$. This case requires a separate consideration.

\end{document}